\documentclass[twoside,a4paper]{article}

\usepackage{a4wide}
\usepackage{fancyhdr}
\usepackage[english]{babel}
\usepackage{amssymb}
\usepackage{amsmath}
\usepackage{theorem}
\usepackage{epsfig}
\usepackage{subfigure}
\usepackage{verbatim}
\usepackage{indentfirst}

\theorembodyfont{\itshape}
\theoremheaderfont{\scshape}
\newtheorem{theor}{Theorem}
\newtheorem{conjc}[theor]{Conjecture}
\newtheorem{propo}{Proposition}[section]
\newtheorem{lemma}[propo]{Lemma}
\newenvironment{proof}{\noindent{\scshape Proof.}}{\hspace{2mm} $\square$ \\}

\newcommand{\R}{\mathbf{R}}
\newcommand{\Z}{\mathbf{Z}}
\newcommand{\super}[1]{^{^{_{\,#1}}}}
\newcommand{\norm}[1]{|\!|#1|\!|}

\DeclareMathOperator{\card}{card}

\DeclareMathOperator{\inter}{int}

\newcounter{picture}

\newcommand{\single}[1]{
\stepcounter{picture}
\vspace{10pt}
\begin{center}
\small{\textbf{Fig. \thepicture \ } \emph{#1}}
\end{center}}

\newcommand{\double}[2]{
\begin{tabular}{cc}
\hspace{190pt} & \hspace{190pt} \\
\stepcounter{picture}
\small{\textbf{Fig. \thepicture \ } \emph{#1}} &
\stepcounter{picture}
\small{\textbf{Fig. \thepicture \ } \emph{#2}}
\end{tabular}}

\newcommand{\mfsi}   {2}
\newcommand{\mfsh}   {3}
\newcommand{\forwd}  {4}
\newcommand{\arrow}  {5}

\begin{document}
\thispagestyle{empty}
\lhead[\thepage]{\textsc{The multitype contact process with frozen states}}
\setcounter{section}{0}
\setcounter{theor}{0}
\setcounter{picture}{0}

\begin{center}
\textbf{\large THE MULTITYPE CONTACT PROCESS WITH FROZEN} \\ \vspace{4pt}
\textbf{\large STATES OR ALLELOPATHY MODEL}
\footnote{\hspace{-16pt} \textit{American Mathematical Society} 1991 \textit{subject classifications}. 60K35} \\
\footnote{\hspace{-16pt} \textit{Key words and phrases}. Allelopathy, interacting particle system, multitype contact process.} \\
\textsc{N. Lanchier} \\ \vspace{16pt}
\end{center}


\begin{abstract}
\noindent The aim of this paper is to study a generalization of the two colors multitype contact process intended to mimic an example of interspecific
 competition called allelopathy.
 Roughly speaking, the strategy of the first species is to inhibit the growth of the second one by freezing the sites it has colonized.
 Our main objective is to prove that this spatial model exhibits some phase transitions comparable with those of the classical multitype contact process.
 In particular, we show that depending on its birth rate each species may conquer the other one.
\end{abstract}


\section{\normalsize\sc\hspace{-10pt}{Introduction}}
\label{mfs-introduction}

\indent Allelopathy denotes a process involving secondary metabolites produced by plants, micro-orga\-nisms, viruses and fungi that influence
 development of biological systems.
 Typically, the strategy of the so-called inhibitory species involved in such a process consists in increasing its competitivity by inhibiting the
 growth of rival species, called susceptible species.
 The phenomenon is now well known and was observed in many ecosystems.
 Oueslati (2003) showed for instance that leaf extracts of \emph{Triticum durum} L. (durum wheat) depresses the germination rate and radicel length
 of \emph{Hordeum vulgare} L. (a barley variety) and \emph{Triticum aestivum} L. (a bread wheat variety).
 Leaves of \emph{Pueraria thunbergiana} also possess allelopathic activity by secreting a substance called xanthoxin that inhibits the root growth of
 \emph{Lepidium sativum} L. (cress) seedlings.
 See Kato-Noguchi (2003).
 Allelopathy cases were also frequently observed in aquatic environment between marine algal species such as \emph{Alexandrium} that has inhibitory
 effects on growth of \emph{Gymnodinium mikimotoi}, \emph{Scrippsiella trochoidea} and \emph{Chaetoceros gracile}.
 See Arzul et al. (1998).
 Usually, these phenomena are characterized by the formation of empty spaces due to the production of toxins that prevent the susceptible species
 from setting up.

\indent The model we introduce to mimic this phenomenon is a continuous time Markov process whose state at time $t$ is a function
 $\xi_t : \Z^d \longrightarrow \{0, 1, 2, 3 \}$.
 A site $x$ in $\Z^d$ is said to be occupied by the \emph{inhibitory species} or \emph{blue particle} if $\xi (x) = 1$, occupied by the
 \emph{susceptible species} or \emph{red particle} if $\xi (x) = 2$, and empty otherwise.
 In the third case, $x$ will be called a \emph{free site} if $\xi (x) = 0$ and a \emph{frozen site} if $\xi (x) = 3$.
 To describe the evolution rules, we now let $\mathcal N \subset \Z^d$ be the set of $y \in \Z^d$ such that $\norm{\,y \,} \leq r$ where $r$ is a
 positive constant and $\norm{\cdot}$ some norm on $\R^d$.
 In other respects, we denote by $f_i \,(x, \xi)$ the fraction of neighbors of $x$ occupied in the configuration $\xi$ by a particle of type $i$,
 where the neighbors of $x$ refer to the translated set $x + \mathcal N$.
 With the notations we have just introduced, we can formulate the transition rates as follows.
 $$ \begin{array}{@{\vspace{2pt}}l@{\qquad}l@{\qquad\qquad\qquad}l@{\qquad}l@{\vspace{2pt}}}
  0 \ \longmapsto \ 1 & \lambda_1 \,f_1 \,(x, \xi) & 1 \ \longmapsto \ 3 & 1 \\
  3 \ \longmapsto \ 1 & \lambda_1 \,f_1 \,(x, \xi) & 3 \ \longmapsto \ 0 & \gamma \\
  0 \ \longmapsto \ 2 & \lambda_2 \,f_2 \,(x, \xi) & 2 \ \longmapsto \ 0 & 1 \end {array} $$
 In particular, the model can be seen as a generalization of the two colors multitype contact process (see Neuhauser (1992)) in which the strategy of
 the blue particles is to inhibit the spread of the red ones by freezing the sites they have colonized for an exponentially distributed amount of time.
 Reciprocally, the multitype contact process can be considered as the extreme case $\gamma = \infty$, i.e. the transition from 3 to 0 is instantaneous.

\indent We now formulate our results and construct step by step the phase diagram of the process.
 First of all, if we start the evolution from the measure $\delta_2$ that concentrates on the all 2's configuration, the process is called the basic
 contact process with parameter $\lambda_2$.
 In such a case, we know that there exists a critical value $\lambda_c \in (0, \infty)$ such that the process converges in distribution to the all
 empty state if $\lambda_2 \leq \lambda_c$, and to a measure $\mu_2$ that concentrates on configurations with infinitely many 2's otherwise.
 See, e.g., Liggett (1999).
 If we start from the all 1's configuration, we have the same result.
 That is, if $\lambda_1 \leq \lambda_c$ then $\xi_t \Rightarrow \delta_0$ while if $\lambda_1 > \lambda_c$ then $\xi_t$ converges to the stationary
 measure $\nu_1 = \lim_{\,t \to \infty} \delta_1 \,S_t$, where $S_t$ denotes the semigroup of the process.
 Moreover, since 3's do not disturb 1's, it is not difficult to prove that $\nu_1 (\xi (x) = 1) = \mu_1 (\xi (x) = 1)$, where $\mu_1$ is the upper
 invariant measure of the basic contact process.

\begin{figure}[ht]
\centering
\mbox{\subfigure{\epsfig{figure = 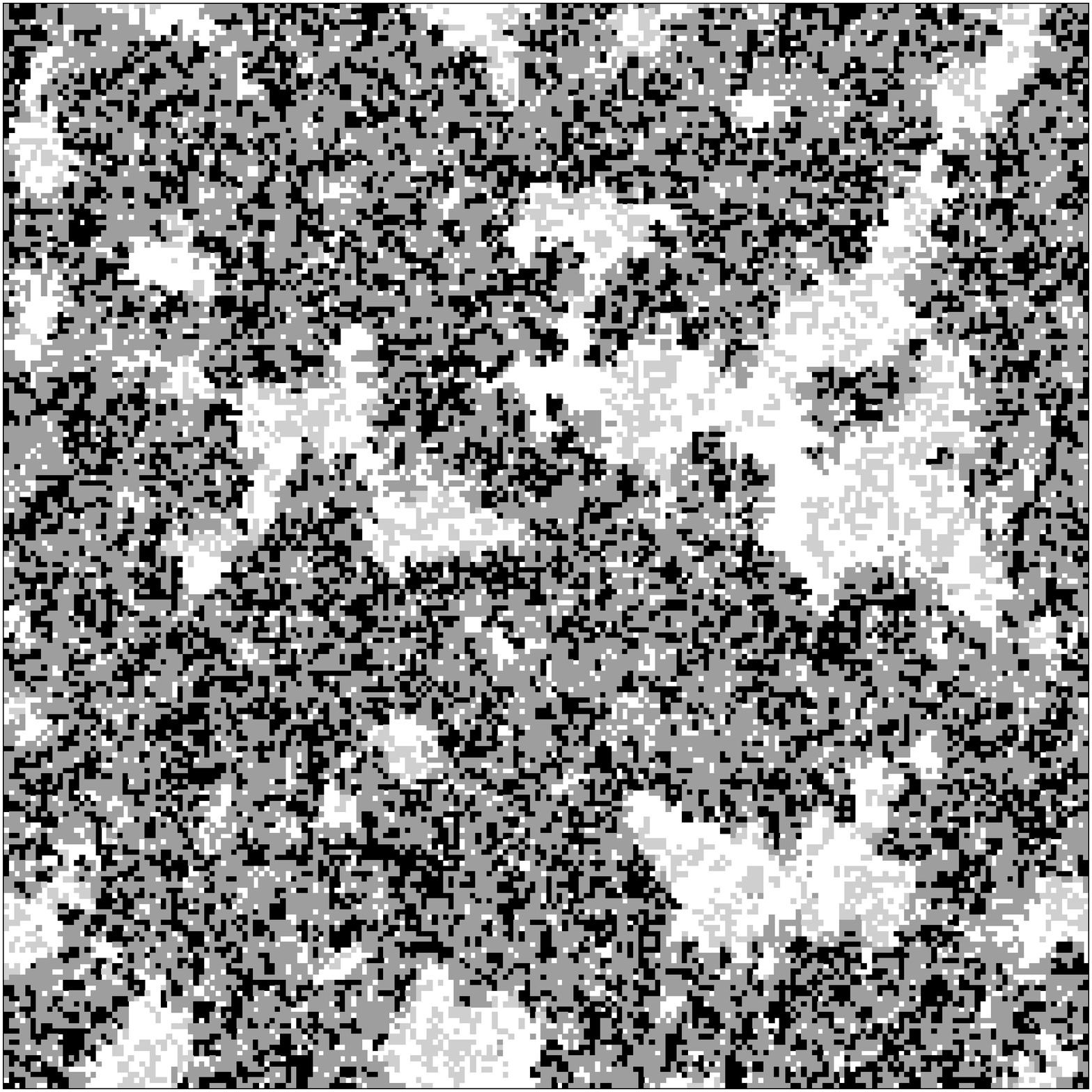, width = 190pt, height = 190pt}} \hspace{20pt}
      \subfigure{\epsfig{figure = 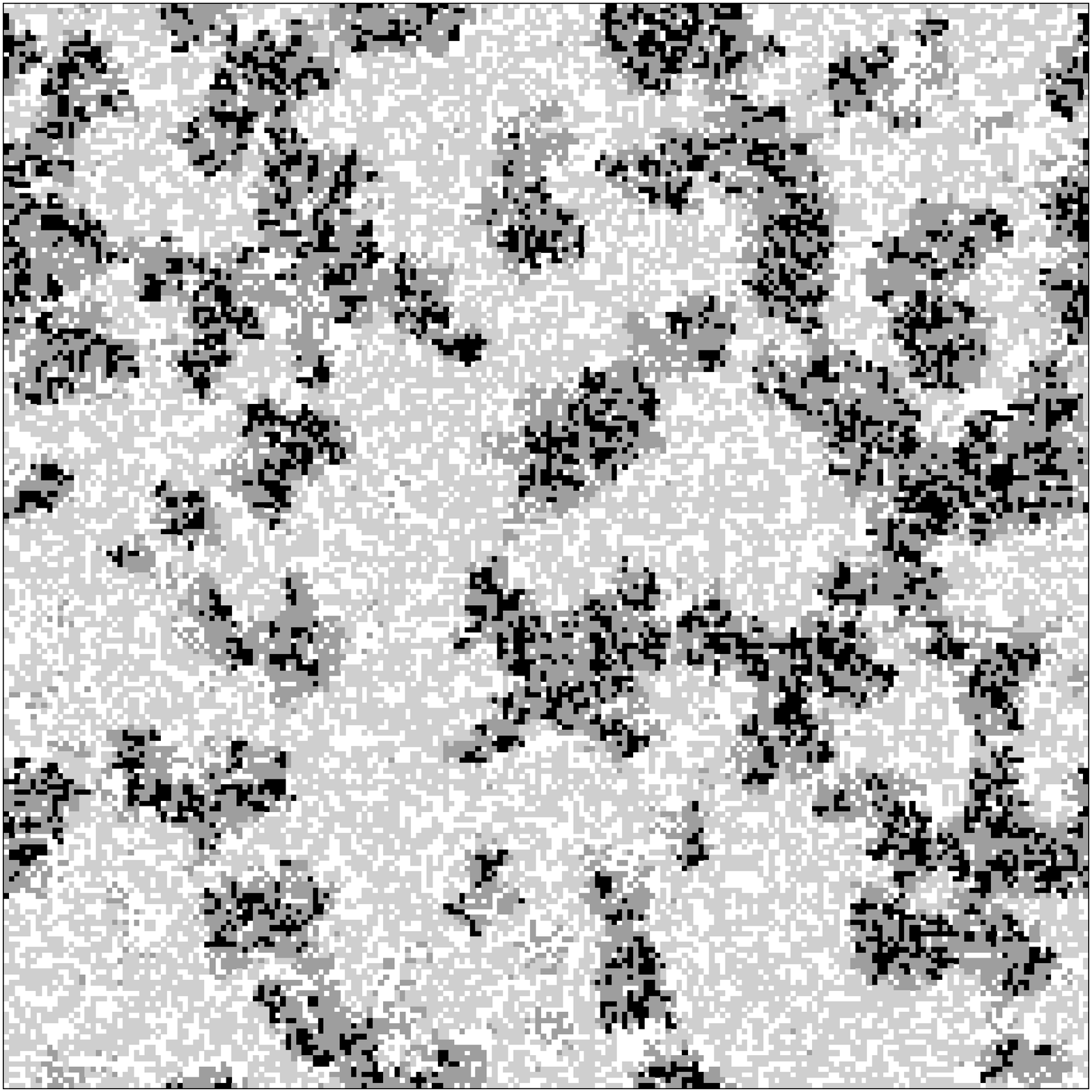, width = 190pt, height = 190pt}}} \\
\stepcounter{picture}
\vspace{10pt}
\parbox{400pt}{\small{\textbf{Fig. \thepicture \ }
\emph{Configurations of the nearest neighbor multitype contact process with frozen state at time $t = 50$ on the square $\{0, \,1, \,\cdots, 199\}^2$.
 Left picture: $\lambda_1 = \lambda_2 = 1,96$ and $\gamma = 0,05$.
 Right picture:  $\lambda_1 = 1,96$, $\lambda_2 = 2,88$ and $\gamma = 0,05$.
 The black points refer to the inhibitory species, the pale gray ones to the susceptible species, and the dark gray ones to frozen sites.}}}
\end{figure}

\noindent To study the competition between both species, we now start the evolution with infinitely many 1's and 2's and suppose that
 $\lambda_1 > \lambda_c$ and $\lambda_2 > \lambda_c$.
 We first set $\gamma_1 < \gamma_2$ and denote by $\xi_t \super{i}$ the allelopathy model with parameters $\lambda_1$, $\lambda_2$ and $\gamma_i$.
 Then, if we think of the processes as being generated by Harris' graphical representation, we may run $\xi_t \super{1}$ and $\xi_t \super{2}$
 simultaneously starting from the same initial configuration in such a way that $\xi_t \super{1}$ has more 1's and fewer 2's than $\xi_t \super{2}$,
 i.e. for any $x \in \Z^d$, if $\xi_t \super{1} (x) = 2$ then $\xi_t \super{2} (x) = 2$, and if $\xi_t \super{2} (x) = 1$ then
 $\xi_t \super{1} (x) = 1$.
 The same argument implies that the process is also monotonous with respect to each of the parameters $\lambda_1$ and $\lambda_2$.

\begin{theor}
\label{monotonous}
 We set $\Theta_t \super{i} = \{\,x \in \Z^d \,;\,\xi_t (x) = i \,\}$ the sites occupied at time $t$ by a type $i$ particle.
 Then, for $i = 1, 2$, the survival probabilities $P \,(\Theta_t \super{i} \neq \varnothing \textrm{ for all } t \geq 0 )$ are monotonous with
 respect to each of the parameters $\lambda_1$, $\lambda_2$ and $\gamma$.
\end{theor}

\noindent In particular, by comparison with the multitype contact process ($\gamma = \infty$), it looks clear that the 1's win as soon as
 $\lambda_1 > \lambda_2$, i.e. $\xi_t \Rightarrow \nu_1$.
 See Theorem 1 in Neuhauser (1992).
 More precisely, we have the following result.

\begin{theor}
\label{mfs-blue}
 We assume that $\xi_0$ is translation invariant.
 If $\lambda_1 > \lambda_2$ or $\lambda_1 = \lambda_2$ and $d \geq 3$, then $\xi_t \Rightarrow \nu_1$ the limit starting from the all 1's configuration.
\end{theor}

\noindent The proof of Theorem \ref{mfs-blue} for $\lambda_1 > \lambda_2$ simply relies on a coupling argument.
 On the other hand, if we now assume that $\lambda_1 = \lambda_2$ and $d \geq 3$, comparison with the multitype contact process just implies
 survival of 1's for any $\gamma > 0$.
 See Neuhauser (1992), Theorem 3.
 In this case, the key of the proof is duality.
 The ancestry of a given point $(x, t)$ in $\Z^d \times \R^+$ exhibits a tree structure that can be divided, as for the ecological succession model
 (see Lanchier (2003)), into several layers, due to the transition $1 \to 3$.
 The idea is that the second layer provides enough 1's to freeze the path of the first ancestor at infinitely many points, which blocks 2's from
 determining the color of $(x, t)$.
 Our approach is roughly the same as in Lanchier (2003), Sect. 4 and Sect. 5, and we think that it may be extended to a larger class of spatial models.

\indent If we now focus on the case $\lambda_1 < \lambda_2$, it is not clear that the 2's win since the dynamics is favorable to the 1's.
 Nevertheless, by refering again to the extreme case $\gamma = \infty$, we can prove that for fixed $\lambda_1 < \lambda_2$ the 2's win in $d = 2$ and for
 $\gamma$ sufficiently large.

\begin{theor}
\label{mfs-red}
 We suppose that $d = 2$ and $\lambda_1 < \lambda_2$.
 Then, there exists a critical value $\gamma_c \in (0, \infty)$ such that $\xi_t \Rightarrow \mu_2$ for any $\gamma > \gamma_c$.
\end{theor}

\noindent The proof of Theorem \ref{mfs-red} is quite simple.
 It essentially relies on the rescaling argument introduced in Durrett and Neuhauser (1997), Sect. 3.
 The idea is to prove that the argument resists to small perturbations of the process.
 That is, if $\gamma > \gamma_c$ then, with probability close to 1, the allelopathy model behaves like the multitype contact process inside a given
 bounded space-time box. \\

\indent To draw the phase diagram of the process, we now define the critical value $\lambda_2 \,(\gamma, \lambda_1)$ as the infimum of
$\lambda_2 \geq 0$ such that the 1's die out, with the convention $\inf \varnothing = \infty$.
 A straightforward application of our theorems then implies that $\lambda_2 \,(\gamma, \lambda_1) \downarrow \lambda_1$ as $\gamma \uparrow \infty$ and
 $\lambda_2 \,(\gamma, \cdot \,)$ is increasing on $\R^+$.
 In conclusion, the picture we finally obtain is given by Fig. \mfsi.
 Unfortunately, we don't know if the 2's may coexist with the 1's for $\lambda_1 < \lambda_2 < \lambda_2 \,(\gamma, \lambda_1)$ (gray part of the picture),
 and if they can win for any $\lambda_1 > \lambda_c$, that is if $\lambda_2 \,(\gamma, \lambda_1) < \infty$ for any $\lambda_1 > \lambda_c$.
 To have an idea of the answers for these questions, we can resort to the mean field theory.
 More precisely, if we pretend that the functions $u_i (x, t) = P \,(\xi_t (x) = i)$ do not depend on the site $x \in \Z^d$, the evolution of the process
 can be formulated in terms of the following ordinary differential equations.
\begin{eqnarray*}
 u_0' & = & u_2 \ + \ \gamma \,u_3 \ - \ \lambda_1 \,u_0 \,u_1 \ - \ \lambda_2 \,u_0 \,u_2 \\
 u_1' & = & \lambda_1 \,u_0 \,u_1 \ + \ \lambda_1 \,u_0 \,u_3 \ - \ u_1 \\
 u_2' & = & \lambda_2 \,u_0 \,u_2 \ - \ u_2 \\
 u_3' & = & u_1 \ - \ \lambda_1 \,u_1 \,u_3 \ - \ \gamma \,u_3
\end{eqnarray*}
 Let $\Omega = \{\,u \,;\,u_i \,\geq \,0, \ u_0 \,+ \,u_1 \,+ \,u_2 \,+ \,u_3 \,= \,1 \}$ be the collection of values we are interested in.
 For fixed $\gamma > 0$, we also denote by $\omega_1$ (resp. $\omega_2$) the set of parameters $\lambda_1$ and $\lambda_2$ such that $\lambda_1 > 1$ and
 $\gamma \,\lambda_2 < \lambda_1 \,(\lambda_1 \,+ \,\gamma \,- \,1)$ (resp. $\lambda_2 > 1$ and $\lambda_2 > \lambda_1$).
 First of all, a straightforward calculation shows that the ODE has a nontrivial fixed point $\bar u$ on the boundary $\bar u_2 = 0$ if and only if
 $\lambda_1 > 1$, where nontrivial means $\neq (1, 0, 0, 0)$.
 Moreover, by studying the eigenvalues of the linearization at point $\bar u$ of the ODE, we can prove that the equilibrium $\bar u$ is attracting if
 $(\lambda_1, \lambda_2) \in \omega_1$, and unstable otherwise, that is the linearization has an unstable direction that points into $\inter (\Omega)$.
 In the same way, if $\lambda_2 > 1$ there exists a nontrivial equilibrium $\bar v$ on the boundary $\bar v_1 = \bar v_3 = 0$ that is attracting if
 $(\lambda_1, \lambda_2) \in \omega_2$, and unstable otherwise.
 Finally, one can prove that the ODE has a fixed point belonging to $\inter (\Omega)$ if and only if $(\lambda_1, \lambda_2) \in \omega_1 \cap \omega_2$.
 In other words, coexistence is possible for the mean field model as soon as
 $\lambda_1 < \lambda_2 < \gamma^{-1} \,\lambda_1 \,(\lambda_1 \,+ \,\gamma \,- \,1)$.

\indent Typically, if the range of the interaction $r$ or the dimension $d$ are sufficiently large, one expects that the process has the same
 properties as the mean field model.
 In particular, relying on the instability of $\bar u$ for $(\lambda_1, \lambda_2) \in \omega_1$, we believe that for any $\lambda_1 > \lambda_c$
 the 2's win for $\lambda_2$ large.
 However, in view of the properties of the mutitype contact process (See Neuhauser (1992)), we don't think that coexistence may occur in an open set
 of the parameters $\lambda_1$ and $\lambda_2$.
 In conclusion, we formulate the following conjecture.

\begin{conjc}
 If the range of the interaction $r \geq r_0$ or the dimension $d \geq d_0$, then for any $\lambda_1$ there exists $\lambda_2 \,(\gamma, \lambda_1) < \infty$
 such that the 1's win if $\lambda_2 < \lambda_2 \,(\gamma, \lambda_1)$ and the 2's win if $\lambda_2 > \lambda_2 \,(\gamma, \lambda_1)$.
\end{conjc}

\begin{figure}[ht]
\centering
\scalebox{0.33}{\input{mfsi.pstex_t}}
\single{Phase diagram.}
\end{figure}


\section{\normalsize\sc\hspace{-10pt}{Graphical representation and duality}}
\label{mfs-construction}

\indent We begin by constructing the process from a collection of Poisson processes in the case $\lambda_1 \geq \lambda_2$.
 For $x$, $y \in \Z^d$, $x - y \in \mathcal N$, let $\{ T_n^{x,y} \,;\, n \geq 1 \}$, $\{ U_n^x \,;\, n \geq 1 \}$ and
 $\{ V_n^x \,;\, n \geq 1 \}$ be the arrival times of Poisson processes with rates $\lambda_1 / \card \mathcal N$, 1 and $\gamma$
 respectively.
 At times $T_n^{x,y}$, we draw an arrow from $x$ to $y$, toss a coin with success probability $(\lambda_1 - \lambda_2) / \lambda_1$, and, if there is a
 success, label the arrow with a 1.
 If $x$ is occupied by a red particle, that the arrow is unlabelled and that $y$ is free, the particle will give birth through this arrow.
 If $x$ is occupied by a blue particle and that $y$ is free or frozen, the site $y$ will be painted in blue.
 Now, at times $U_n^x$, we put a cross $\times$ at $x$.
 The effect of a $\times$ is to kill both species, i.e. a red particle gives way to a free site and a blue one to a frozen site.
 Finally, at times $V_n^x$, we put a dot $\bullet$ at $x$ to indicate that a frozen site becomes free.
 A result of Harris (1972) implies that such a graphical representation can be used to construct the process starting from any initial configuration
 $\xi_0 \in \{0, 1, 2, 3 \}^{\Z^d}$.
 For an example of realization of the process, see Fig. \mfsh.

\begin{figure}[ht]
\centering
\scalebox{0.35}{\input{mfsh.pstex_t}} \\
\stepcounter{picture}
\vspace{20pt}
\parbox{370pt}{\small{\textbf{Fig. \thepicture \ }
\emph{Harris' graphical representation. The black lines refer to blue particles, the pale gray ones to red particles,
 the dark gray ones to frozen sites, and the dotted ones to free sites.}}}
\end{figure}

\indent After constructing the graphical representation, we can now define the dual process.
 We say that two points $(x, s)$ and $(y, t)$ in $\Z^d \times \R^+$ are \emph{connected} or that there is a \emph{path} from $(x, s)$ to $(y, t)$
 if there exists a sequence of times $s_0 = s < s_1 < s_2 < \cdots < s_n < s_{n + 1} = t$ and spatial locations
 $x_0 \, = \,x, \ x_1, \ x_2, \ \cdots \ x_n \, = \,y$ so that
\begin{enumerate}
\item [(1)] for $1 \leq i \leq n$, there is an arrow from $x_{i - 1}$ to $x_i$ at time $s_i$ and
\item [(2)] for $0 \leq i \leq n$, the vertical segments $\{x_i\} \times (s_i, s_{i + 1})$ do not contain any $\times$'s.
\end{enumerate}
 If (1) only is satisfied, we say that there is a \emph{weak path} from $(x, s)$ to $(y, t)$.
 Finally, we say that there exists a \emph{dual path} from $(x, t)$ to $(y, t - s)$, $0 \leq s \leq t$, if there is a path from $(y, t - s)$ to $(x, t)$.
 In other words, dual paths move against the direction of time and arrows.
 We then define the \emph{dual process} by setting
 $$ \hat \xi_s\super{(x, t)} \ = \ \{\, y \in \Z^d \,;\, \textrm{there is a dual path from $(x, t)$ to $(y, t - s)$} \} $$
 for any $0 \leq s \leq t$.
 First of all, we can observe that $\{(\hat \xi_s\super{(x, t)}, s) \,;\, 0 \leq s \leq t \,\}$, which is the set of points in $\Z^d \times [0, t \,]$
 that are connected with $(x, t)$, exhibits a tree structure.
 As for the mutitype contact process, such a structure allows to equip the dual process $\hat \xi_s\super{(x, t)}$ with an ordered relation for which
 the members are arranged according to the order they determine the color of $(x, t)$.
 See e.g., Neuhauser (1992), Sect. 2.
 From now on, the tree
 $$ \Gamma \ = \ \{(\hat \xi_s\super{(x, t)}, s) \,;\,0 \leq s \leq t \,\} $$
 will be called the \emph{upper tree starting at} $(x, t)$ and the elements of $\hat \xi_s\super{(x, t)}$ the \emph{upper ancestors}.
 We will denote by $\hat \xi_s\super{(x, t)} (n)$ the $n$-th member of the ordered ancestor set.
 For an example of ancestor hierarchy, see Fig. \forwd.
 If there is a weak path from $(y, r)$ to $(x, t)$, the tree starting at $(y, r)$ will be called a \emph{lower tree} and the elements of
 $\hat \xi_s\super{(y, r)}$ the \emph{lower ancestors}.
 Finally, the first upper ancestor, that is $\hat \xi_s\super{(x, t)} (1)$, will be called the \emph{distinguished particle}.

\indent To conclude this section, we now describe in greater detail a general method to determine the color of $(x, t)$ in the case
 $\lambda_1 \geq \lambda_2$.
 Contrary to the multitype contact process, the state of some sites (free or frozen) depends on the structure of the lower trees.
 So, for more convenience, we first assume that this information is known and say that an arrow is \emph{forbidden for the red} if its target site is frozen
 or if it is a 1-arrow.
 We will see further that our results can be proved by focusing only on the first two layers of the tree structure.
 First of all, if the distinguished particle lands on a blue site, it will paint $(x, t)$ in blue.
 If it lands on a red site and does not cross any arrow forbidden for the red, it will paint $(x, t)$ in red.
 On the other hand, if it lands on a red site and crosses at least one arrow forbidden for the red, we follow the path it takes to reach $(x, t)$ until
 we first meet a forbidden arrow.
 Then, we look backwards in time starting from the location where this arrow is attached and discard all the ancestors of this point, that is the next few
 members of the ancestor hierarchy.
 Finally, we start over again with the next remaining ancestor, and so on.
 In the same way, if the distinguished particle lands on a free or frozen site, we start afresh with the second ancestor.

\indent For instance, in Fig. \forwd \ the distinguished particle lands on a red site and crosses an arrow that points at the frozen site $x - 3$ so
 fails in painting $(x, t)$ any color.
 The same holds for the second ancestor that takes a 1-arrow to reach $(x, t)$.
 In other respects, the third ancestor lives at time 0 on the tree starting at the location where this 1-arrow is attached, so we look at the
 fourth ancestor.
 This last one finally succeeds in painting $(x, t)$ in blue.


\section{\normalsize\sc\hspace{-10pt}{Proof of Theorem \ref{mfs-blue}}}
\label{mfs-duality}

\indent In this section, we will prove Theorem \ref{mfs-blue} beginning with the case $\lambda_1 > \lambda_2$.
 The strategy of the proof is quite simple.
 It relies on a basic coupling between the allelopathy model and the multitype contact process.
 First of all, recall that the multitype contact process can be considered as the extreme case $\gamma = \infty$, i.e. the transition from 3 to 0 is
 instantaneous.
 In that case, it is known that if $\xi_0 \in \{0, 1, 2 \}^{\Z^d}$ is translation invariant and $\lambda_1 > \lambda_2$ then starting from infinitely
 many 1's and 2's the first species win the competition, i.e. $\xi_t \Rightarrow \mu_1$ the upper invariant measure of the basic contact process.
 See Theorem 1 in Neuhauser (1992).
 Now the intuitive idea is that the smaller $\gamma$ is the more the blue particles are competitive.
 To specify this, we couple both processes together considering two identical copies $G_1$ and $G_2$ of the graphical representation.
 In the first copy, the particles evolve according to the rules given above.
 In the second one, the $\times$'s permute both species into free site, whereas the dots have no more effect on the particles so $G_2$ can be seen as the
 graphical representation of the multitype contact process.
 Then, we may run both processes $\xi_t\super{1}$ and $\xi_t\super{2}$ simultaneously on $G_1$ and $G_2$ starting from the same configuration
 in such a way that $\xi_t\super{1}$ has more 1's and fewer 2's than $\xi_t\super{2}$, i.e. for any point $(x, t)$ in $\Z^d \times \R^+$, if
 $\xi_t\super{1} (x) = 2$ then $\xi_t\super{2} (x) = 2$, and if $\xi_t\super{2} (x) = 1$ then $\xi_t\super{1} (x) = 1$.
 This together with Theorem 1 of Neuhauser (1992) implies the first part of Theorem \ref{mfs-blue}. \\

\begin{figure}[ht]
\centering
\scalebox{0.4}{\input{forwd.pstex_t}} \hspace{50pt}
\scalebox{0.4}{\input{arrow.pstex_t}}
\double{\hspace{25pt}}{\hspace{-25pt}}
\end{figure}

\indent The extension of this result for $\lambda_1 = \lambda_2$ and $d \geq 3$ is more difficult.
 In this case, the strategy consists in proving that if the upper tree $\Gamma$ lives forever, the distinguished particle jumps infinitely often
 on a frozen site.
 To do this, we will focus on the second layer of the tree structure to show that, conditioned on survival of the upper tree, the number of blue
 particles that freeze the sites visited by the first ancestor can be made arbitrarily large.
 We will finally conclude by exhibiting an upper ancestor that will paint $(x, t)$ in blue.

\indent From now on, we denote by $\lambda$ the common value of $\lambda_1$ and $\lambda_2$ and suppose that $\Gamma$ lives forever.
 We can observe that such an event obviously occurs with positive probability since $\lambda$ is supercritical.
 For more convenience, we reverse the arrows and time by letting $\tilde s = t - s$.
 The main objective is to prove that the number of frozen sites visited by the first ancestor tends to infinity as $t \uparrow \infty$.
 To specify this result, we start by defining more precisely our framework.
 First of all, we follow the path of the distinguished particle starting from $(x, \tilde 0)$ and, for any $n \geq 1$, denote by
 $\beta_n = \beta_n (x, t)$ the $n$-th crossed arrow.
 See Fig. \arrow \ for a picture.
 Now let $z_n = z_n (x, t)$ and $\tilde s_n = \tilde s_n (x, t)$ be respectively the starting site and the temporal location of the arrow
 $\beta_n (x, t)$ and denote by $N (x, t)$ the number of arrows $\beta_n (x, t)$ that start at a frozen site, i.e.
 $$ N (x, t) \ = \ \card \,\{ n \geq 1 \,;\, \xi_{\tilde s_n} (z_n) = 3 \}. $$
 By construction of the sequence $\beta_n (x, t)$, it is clear that $N (x, t)$ also denotes the number of frozen sites visited by the distinguished
 particle.
 In particular, the main result we have to prove can be formulated as follows.

\begin{propo}
\label{mfs-frozen}
 If $\lambda_1 = \lambda_2$ and $d \geq 3$ then \ $\lim_{\,t \to \infty} N (x, t) = \infty$ \ a.s.
\end{propo}

\noindent The intuitive idea of the proof is that the second layer of the tree structure provides enough 1's to freeze the
 path of the distinguished particle at infinitely many points so that we can focus only on this one.
 To specify this point, we proceed in two steps.

\indent First of all, we denote by $\tilde \sigma_n = \tilde \sigma_n (x, t)$ the arrival time of the first $\times$ located under the starting point of
 $\beta_n (x, t)$, i.e.
 $$ \tilde \sigma_n (x, t) \ = \ \min_{k \geq 1} \,\{\,\tilde U_k^{z_n} \textrm{ such that } \tilde U_k^{z_n} \geq \tilde s_n (x, t)\} $$
 and let $\Gamma_n = \Gamma_n (x, t)$ be the lower tree starting at $(z_n, \tilde \sigma_n)$, that is
 $$ \Gamma_n \ = \ \{(y, \tilde s) \in \Z^d \times [\,\tilde \sigma_n, \infty) \,;\, \textrm{there is a path from
 $(z_n, \tilde \sigma_n)$ to $(y, \tilde s)$}\}. $$
 See Fig. \arrow \ for a picture.
 We say that the $n$-th lower tree is \emph{favorable to the 1's} if the following two conditions are satisfied.
\begin{enumerate}
\item [(1)] $\Gamma_n (x, t)$ lives forever and
\item [(2)] the vertical segment $\{z_n\} \times (\tilde s_n, \tilde \sigma_n)$ does not contain any dots.
\end{enumerate}
 As we will see further, the properties (1) and (2) will give us a good opportunity to freeze the site $z_n$ at time $\tilde s_n$ so the first step
 of the proof is to show that there exist infinitely many lower trees $\Gamma_n$ that are favorable to the 1's.

\begin{lemma}
\label{mfs-good}
 Let $G_n$ be the event that the lower tree $\Gamma_n$ is favorable to the 1's.
 If the upper tree lives forever then $P \,(\,\limsup_{n \to \infty} G_n \,) = 1$.
\end{lemma}

\begin{proof}
 To begin with, denote by $A_n$ the event that $\Gamma_n$ lives forever and by $B_n$ the event that $\{z_n\} \times (\tilde s_n, \tilde \sigma_n)$
 does not contain any dots.
 The first step is to prove that for any $n \geq 1$, there exists a.s. an integer $m \geq n$ such that $A_m$ occurs.
 To do this, we set $\Gamma_{n_1} = \Gamma_n$ and while $\Gamma_{n_k}$ is bounded we denote by $\Gamma_{n_{k + 1}}$ the first lower tree that
 borns after $\Gamma_{n_k}$ dies.
 Note that if $A_{n_k}$ does not occur then $\Gamma_{n_{k + 1}}$ is well defined and the event $A_{n_{k + 1}}$ is determined by parts of the graph
 that are after $\Gamma_{n_k}$ dies so that $A_{n_k}$ and $A_{n_{k + 1}}$ are independent.
 More generally, since the trees $\Gamma_{n_1}$, $\Gamma_{n_2}$, $\cdots$, $\Gamma_{n_{k + 1}}$ are disjoint, the events $A_{n_1}$, $A_{n_2}$,
 $\cdots$, $A_{n_{k + 1}}$ are independent.
 Moreover, the probability that $A_{n_k}$ occurs is given by the survival probability $p_{\lambda}$ of the basic contact process with parameter $\lambda$
 starting from one infected site so that
\begin{eqnarray*}
 P \,(\,A_n^c \,\cap \,A_{n + 1}^c \,\cap \,\cdots \,) & \leq &
   \lim_{k \to \infty} \ P \,(\,A_{n_1}^c \,\cap \,A_{n_2}^c \,\cap \,\cdots \,\cap \,A_{n_k}^c \,) \\ & \leq &
   \prod_{k = 1}^{\infty} \,P \,(A_{n_k}^c) \ = \ \lim_{k \to \infty} \ (1 - p_{\lambda})^k \ = \ 0
\end{eqnarray*}
 as soon as $\lambda > \lambda_c$.
 In particular,
 $$ P \,(\,\limsup_{n \to \infty} A_n \,) \ = \ \lim_{n \to \infty} \ P \,(\,A_n \,\cup \,A_{n + 1} \,\cup \,\cdots \,) \ = \ 1. $$
 This proves that, with probability 1, there exist infinitely many lower trees that live forever.
 In other respects, since $\tilde \sigma_n - \tilde s_n$ is exponentially distributed with parameter $\lambda$, we can state that for any $n \geq 1$
 $$ P \,(B_n) \ = \ P \,(\,\tilde \sigma_n - \tilde s_n \leq V_1^{z_n} \,) \ = \ \lambda \,\gamma^{-1} \,(\lambda + \gamma)^{-1} \ > \ 0. $$
 By independence, we can finally conclude that $P \,(\,\limsup_{n \to \infty} A_n \,\cap \,B_n \,) = 1$.
\end{proof}

\indent To complete the proof of Proposition \ref{mfs-frozen}, we now consider for any $n \geq 1$ and $\tilde s \geq \tilde \sigma_n$ the time
 translation dual process
 $$ \hat \xi_s \super{(z_n, \tilde \sigma_n)} \ = \ \{\, y \in \Z^d \,;\, \textrm{there is a path from $(z_n, \tilde \sigma_n)$ to $(y, \tilde s)$} \} $$
 and denote by $\gamma_s (n)$ the associated distinguished particle, that is the first ancestor of $(z_n, \tilde \sigma_n)$.
 First of all, we can observe that if the lower tree $\Gamma_n$ lives forever then $\gamma_s (n)$ is well defined for any $\tilde s \geq \tilde \sigma_n$.
 Moreover, if we suppose that $\Gamma_n$ is favorable to the 1's and that $\gamma_s (n)$ lands at time $\tilde t$ on a blue site it is clear,
 in view of the condition (2) above, that $z_n$ will be frozen at time $\tilde s_n$.
 In particular, letting $\Gamma_{n_k}$ be a subsequence of favorable trees given by Lemma \ref{mfs-good}, the proof of Proposition \ref{mfs-frozen} can
 be completed with the following lemma.

\begin{lemma}
 Let $\tilde \Omega_s = \{\gamma_s (n_k) \,;\,k \geq 1\}$ and $\tilde \Theta_s \super{1}$ be the set of sites occupied at time $\tilde s$ by a 1.
 If $\xi_0$ is translation invariant and $d \geq 3$ then, starting from infinitely many 1's,
 $$ \lim_{t \to \infty} \ \card \,(\tilde \Omega_t \,\cap \,\tilde \Theta_t \super{1}) \ = \ \infty \quad \textrm{a.s.} $$
\end{lemma}

\begin{proof}
 The key of the proof is transience.
 More precisely, a straightforward application of Lemma 5.4 of Lanchier (2003) implies that $\lim_{\,t \to \infty} \,\card \,(\tilde \Omega_t) = \infty$
 a.s. as soon as $d \geq 3$.
 We then conclude with Lemma 9.14 of Harris (1976) and translation invariance of the initial configuration.
\end{proof}

\indent To conclude the proof of Theorem \ref{mfs-blue}, we now restore the direction of time and arrows and construct a sequence of upper
 ancestors $\zeta_t\super{(x, t)} (k)$, $k \geq 0$, that are candidates for painting $(x, t)$ in blue.
 The first member, that is $\zeta_t\super{(x, t)} (0)$, is the distinguished particle $\hat \xi_t\super{(x, t)} (1)$.
 Next, we rename the sequence of frozen points $(z_k, s_k)$, $k \geq 1$, visited by the distinguished particle by going forward in time.
 For fixed $k \geq 1$, we now look backwards in time starting from the location where the arrow $\beta_k (x, t)$ is attached and discard all the
 offspring of this particular point.
 We then define $\zeta_t \super{(x, t)} (k)$ as the first upper ancestor that is left after discarding.
 We now set $\zeta_t = \{ \zeta_t\super{(x, t)} (k) \,;\, k \geq 0 \}$.
 By Proposition \ref{mfs-frozen} and the fact that the upper tree $\Gamma$ is linearly growing in time, the cardinality of $\zeta_t$ can be made
 arbitrarily large by choosing $t$ large enough.
 In particular, a new application of Lemma 9.14 in Harris (1976) gives us that
 $$ \lim_{t \to \infty} \ P \,(\,\zeta_{t - 1} \,\cap \,\Theta_1 \super{1} = \,\varnothing \,) \ = \ 0. $$
 where $\Theta_s \super{1}$ denotes the set of sites occupied at time $s$ by a 1.
 Hence, there exists at least one candidate that lands on a blue site.
 Denote by $\zeta_t \super{(x, t)} (k_0)$ the first one in the hierarchy.
 Since by construction the arrow $\beta_{k_0} (x, t)$ is forbidden for the red particles (we recall that at time $s_{k_0}$ the site $z_{k_0}$ is frozen),
 the upper ancestor $\zeta_t \super{(x, t)} (k_0)$ will finally paint $(z_{k_0}, s_{k_0})$ and so $(x, t)$ in blue.
 This completes the proof of Theorem \ref{mfs-blue}.


\section{\normalsize\sc\hspace{-10pt}{Proof of Theorem \ref{mfs-red}}}
\label{mfs-rescaling}

\indent In this section, we assume that $d = 2$, set $\lambda_1 < \lambda_2$ and prove that there exists $\gamma_c > 0$ such that for any
 $\gamma > \gamma_c$ the 2's win.
 In view of the dynamics, the survival of 2's for $\lambda_1 < \lambda_2$ is not clear and tools as basic coupling and duality techniques fail in proving
 Theorem \ref{mfs-red}.
 We will first rely on the rescaling argument in Durrett and Neuhauser (1997), Sect. 3, applied to the multitype contact process, and then prove that taking
 $\gamma$ large enough does not affect too much the process.
 We start by introducing the suitable space and time scales.
 We let $L$ define a positive integer and, for any $z = (z_1, z_2)$ in $\Z^2$, set
 $$ \Phi (z) \ = \ (L z_1, L z_2), \quad \qquad B \ = \ [- L, L \,]^2, \quad \qquad B (z) \ = \ \Phi (z) + B. $$
 Moreover, we tile $B (z)$ with $L \super{0,1} \times L \super{0,1}$ squares by setting
 $$ \pi (w) \ = \ (L \super{0,1} w_1, L \super{0,1} w_2),  \hspace{15pt}  D \ = \ (- L \super{0,1} / \,2, L \super{0,1} / \,2 \,]^2, $$
 $$ D (w) \ = \ \pi (w) + D,  \hspace{15pt}  I_z \ = \ \{\,w \in \Z^2 \,;\, D (w) \subseteq B (z) \}. $$
 We say that $B (z)$ is \emph{good} if $B (z)$ is void of 1's and has at least one particle of type 2 in each of the squares $D (w)$ for $w \in I_z$.
 For $z = (z_1, z_2) \in \Z^2$ with $z_1$ and $z_2$ both even for even $k$, and $z_1$ and $z_2$ both odd for odd $k$, we say that $(z, k)$
 is \emph{occupied} if $B (z)$ is good at time $kT$, where $T$ is an integer to be picked later on.
 Moreover, we require this event to occur for the process restricted to the space box $[- ML, ML \,]^2 + \Phi (z)$.
 We start by assuming that $\xi_t$ is the multitype contact process, that is $\gamma = \infty$.

\begin{propo}[Durrett and Neuhauser]
\label{mfs-percolation}
 Assume that $\lambda_2 > \lambda_1$ and $T = L^2$.
 For any $\delta > 0$, there exist large enough $L$ and $M$ such that the set of occupied sites dominates the set of open sites in an $M$-dependent
 oriented percolation process with parameter $1 - \delta$.
\end{propo}

\noindent See Durrett and Neuhauser (1997), Proposition 3.1 and Lemma 3.7.
 To extend the comparison result to $\gamma > 0$ sufficiently large, we just need to prove that, with probability close to 1, the process behaves like
 the multitype contact process (i.e. none of the 2's is blocked by a frozen site) inside the space-time box
 $$ \Phi (z) \ + \ [- M L \,/ \,3, \,M L \,/ \,3 \,]^2 \ \times \ [0, T] $$
 See Lemma 3.7 in Durrett and Neuhauser (1997).
 Clearly, this occurs if we free all the target sites by setting a dot under each tip of arrow.
 So, letting $\iota (x, t)$ be the number of arrows that point at $x$ by time $t$ and $\kappa (z) = \Phi (z) + [- M L \,/ \,3, \,M L \,/ \,3 \,]^2$,
 and decomposing according to whether $\iota (x, T) > 2 \,\lambda_2 T$ or $\iota (x, T) \leq 2 \,\lambda_2 T$, we obtain
\begin{eqnarray*}
 P \,(\,\textrm{any of the 2's is blocked} \,)
  & \leq & \sum_{x \in \kappa (z)} P \,(\iota (x, T) > 2 \,\lambda_2 T) \ + \ 2 \,\lambda_2 T \,\sum_{x \in \kappa (z)} P \,(\,V_1^x < U_1^x) \\
  & \leq & (2 \,/ \,3) \,M L \left(C \,e^{- \beta T} \, + \,2 \,\lambda_2 T \,(\gamma \,(\gamma + 1))^{-1} \right) \ \leq \ \delta
\end{eqnarray*}
 for $T$ and $\gamma$ sufficiently large and appropriate $C < \infty$ and $\beta > 0$.
 Comparison with $M$-dependent oriented site percolation and Proposition \ref{mfs-percolation} then imply Theorem \ref{mfs-red}. \\


\noindent\textsc{Acknoledgement}.
 I would like to thank Claudio Landim, Olivier Benois and Roberto Fern\'andez for their advice, and Pierre Margerie for the biological interpretation.


\vspace{5pt}
\noindent\small{\textsc{UMR 6085, Universit\'e de Rouen, 76128 Mont Saint Aignan, France. \\
 E-mail}: nicolas.lanchier@univ-rouen.fr}

\end{document}